
%
%

\input amstex
\documentstyle{amsppt}
\NoBlackBoxes
\def\card{\operatorname{card}}
\def\rmod{\operatorname{Mod-}}
\def\bmax{\operatorname{bmax}}
\def\amax{\operatorname{amax}}
\def\Zg{\operatorname{Zg}}
\def\Spec{\operatorname{Spec}}

\def\Image{\operatorname{Im}}
\def\End{\operatorname{End}}

\def\ker{\operatorname{Ker}}
\topmatter
\title
Spectra of the $\Gamma$-invariant of uniform modules
\endtitle
\author
Saharon Shelah and Jan Trlifaj
\endauthor

\thanks First author publication number 693.
Research of the second author supported by a Fulbright
Scholarship at UCI. His thanks are due to Professor Paul Eklof for
many stimulating discussions on the subject, and for his constant
help. Thanks are also due to Rutgers University for
supporting the second author's trip to Rutgers. \endthanks

\abstract For a ring $R$, denote by $\Spec _\Gamma(\kappa,R)$
the $\kappa$-spectrum of the $\Gamma$-invariant of strongly uniform
right $R$-modules. Recent realization techniques of Goodearl
and Wehrung show that $\Spec _\Gamma(\aleph_1,R)$ is full
for a suitable von Neumann regular algebra $R$, but the techniques
do not extend to cardinals $\kappa > \aleph_1$. By a direct
construction, we prove that for any field $F$ and any regular
uncountable cardinal $\kappa$ there is an $F$-algebra $R$ such that
$\Spec _\Gamma(\kappa,R)$ is full.
We also derive some consequences for the $\Gamma$-invariant of
strongly dense lattices of two-sided ideals, and for the complexity of
Ziegler spectra of infinite dimensional algebras.
\endabstract
\endtopmatter

\document
The $\Gamma$-invariant method introduced by Eklof in \cite{E1}
and \cite{E2} provides an efficient tool for classification of
algebraic objects which are defined by existence of infinite
filtrations of particular forms.
The method has been used to develop a structure theory of almost free
groups \cite{EM}, uniserial modules \cite{Sa},
and bilinear spaces \cite{A}, \cite{BFS}.

More recently, $\Gamma$-invariants were defined
also in the dual setting, for objects possessing dual
filtrations. This resulted in a classification of dense lattices
\cite{ET}, and of strongly uniform modules \cite{T1},
\cite{T2}.

For a regular uncountable cardinal $\kappa$, denote by $B(\kappa)$
the Boolean algebra consisting of all subsets of $\kappa$ modulo
the filter of subsets containing a closed unbounded set.
The $\Gamma$-invariant of objects of dimension $\kappa$ takes
values in $B(\kappa)$. The value measures a caveat for an
object of dimension $\kappa$ to have a certain algebraic
property. For example, for almost free groups, the property is
``to be a free group'' \cite{E3}. For bilinear spaces, the
property is ``to decompose orthogonally into subspaces of
dimension $< \kappa$'' \cite{BFS}. For dense lattices,
it is ``to be relatively complemented'' \cite{ET}, etc.

For each $\Gamma$-invariant, two natural problems arise:

\smallskip

(1) \, Given a regular uncountable cardinal $\kappa$ and $i \in
B(\kappa)$, is there an object of dimension $\kappa$
whose $\Gamma$-invariant value equals $i$ ?

\smallskip

The set of all $i \in B(\kappa)$ for which the answer to (1) is
positive is called the {\it $\kappa$-spectrum} of the $\Gamma$-invariant,
and denoted by $\Spec _\Gamma(\kappa)$.
The $\kappa$-spectrum is said to be {\it full} provided that
$\Spec _\Gamma(\kappa) = B(\kappa)$, \cite{BFS}.

\smallskip

(2) For $i \in \Spec _\Gamma(\kappa)$, describe all the objects of
dimension $\kappa$ whose $\Gamma$-invariant value equals $i$.

\smallskip

Solutions to problems (1) and (2) depend substantially on the
particular form of the $\Gamma$-invariant.

For almost free groups, the $\kappa$-spectrum is full for each
$\kappa = \aleph_n$, $n < \omega$, \cite{M, Theorem 5.6},
but the fullness for $\kappa = \aleph_{\omega^2 + 1}$ is independent
of ZFC \cite{EM}, \cite{MS}. For bilinear spaces,
the $\kappa$-spectrum is full for $\kappa = \aleph_1$ \cite{A},
but it is not full for any regular $\kappa \geq \aleph_2$ \cite{S, Lemma 2}.
For dense lattices, the $\kappa$-spectrum is full for all regular
uncountable cardinals $\kappa$ \cite{ET, Theorem 1.15}.

Since isomorphic objects have the same value of the
$\Gamma$-invariant, fullness of the $\kappa$-spectrum always
implies that there exist many (at least $2^\kappa$) non-isomorphic
objects of dimension $\kappa$. In that case, (2) gives a strategy
for a fine classification of all objects of dimension $\kappa$.

\medskip

In the present paper, we provide a complete solution to problem (1)
for the $\Gamma$-invariant of strongly uniform modules
introduced in \cite{T1}. Answering the questions
of \cite{T1, \S3, Problem 3}, \cite{ET, \S2} and \cite{T2, \S2},
we prove that the $\kappa$-spectrum is full for all regular
uncountable cardinals $\kappa$. Our main result is as follows:

\smallskip

\proclaim{Theorem 2.7} Let $\lambda$ be an uncountable
cardinal and $F$ be a field. Then there exists an $F$-algebra $R$
such that for any regular uncountable cardinal $\kappa \leq
\lambda$ and any $i \in B(\kappa)$ there is a strongly uniform
module $L \in \rmod R$ such that $\End _R (L) = F$ and
$\Gamma (L) = i$. In particular, $\Spec _\Gamma(\kappa,R)$ is
full. \endproclaim

\smallskip

Section 1 contains basic facts about strongly uniform
modules. The proof of Theorem 2.7 is presented in Section 2.
In Section 3, we deal with consequences for the
$\Gamma$-invariant of two-sided ideal lattices.
We also relate our construction to the Goodearl-Wehrung one
(cf\. \cite{GW, Theorem 4.4} and
\cite{T2, Theorem 2.4}). The latter works only for $\kappa
= \aleph_1$, but provides for additional properties of the algebras
and modules. In Section 4, we derive consequences for the
structure of Ziegler spectra of infinite dimensional algebras.

\bigskip

\heading 1. \, Strongly uniform modules \endheading

\medskip

Let $R$ be an associative ring with unit. Denote by $L_2(R)$ the
lattice of all two-sided ideals of $R$, and by $\rmod R$ the category
of all (unitary right $R$-) modules.
If $M \in \rmod R$, then $\End_R(M)$ denotes the
endomorphism ring of $M$. (Endomorphisms are always
written as acting on the opposite side from scalars).

A non-zero module $U \in \rmod R$ is called {\it uniform}
provided that $V \cap W \neq 0$ for all non-zero submodules $V$
and $W$ of $U$. So uniform modules coincide with non-zero
submodules of indecomposable injective modules. Uniform
modules play an important role in module theory: for example,
they form building blocks for Goldie dimension theory of
modules, \cite{MR}. (For the model-theoretic role of
injective uniform modules, we refer to \cite{P1} and \cite{P2};
see also Section 4.)

A trivial sufficient condition for uniformity of a module over
an arbitrary ring is the existence of a minimal non-zero submodule.
Such uniform modules are called {\it cocyclic}. Cocyclic modules
are exactly the strongly uniform modules of dimension 1 in the
sense of the following

\definition{Definition 1.1} \rm Let $R$ be a ring and $U \in \rmod R$.
A sequence of non-zero submodules of~$U$, $\Cal U = ( U_\alpha \mid
\alpha < \kappa )$, is called a {\it c.d.c.} in $U$ provided that
$\Cal U$ is

\smallskip

- continuous ($U_0 = U$, and $U_\alpha = \cap_{\beta < \alpha}
U_\beta$ for all limit ordinals $\alpha < \kappa$),

- strictly decreasing ($U_{\alpha + 1} \subset U_\alpha$ for all
$\alpha < \kappa$), and

- cofinal (for each non-zero submodule $V \subseteq U$ there is $\alpha
< \kappa$ such that $U_\alpha \subseteq V$).

\smallskip

$U$ is {\it strongly uniform} provided that there is a c.d.c\. in
$U$. The ordinal $\kappa$ is called the {\it length} of $\Cal U$.
The least ordinal $\kappa$ such that there is a c.d.c\.
$\Cal U$ of length $\kappa$ in $U$ is called the {\it dimension} of $U$.
\enddefinition

It is easy to see that any strongly uniform module $U$ is uniform, and
either $d = 1$ or $d$ is a regular infinite cardinal, where $d$ is the
dimension of $U$.

Clearly, $d = 1$ iff $U$ is cocyclic. Moreover, any module with a countable
submodule lattice is uniform iff it is strongly uniform.
This is not true in general: if $R = k[x]$ is the polynomial ring of
one variable $x$ over a field $k$ then $U = R$ is uniform, but $U$
is strongly uniform iff $k$ is countable, cf\. \cite{T1, \S 2}.

\definition{Definition 1.2} Let $U$ be a strongly uniform module.
Let $0 \neq V \subset W \subseteq U$. Then
$W$ is {\it complemented} over $V$ (in $U$) provided that there is
a submodule $X \subseteq U$ such that $W \cap X = V$ and $W + X = U$.
For example, $U$ is complemented over any $0 \neq V \subset U$.
\enddefinition

\medskip

Also the case of the least infinite dimension,
$d = \aleph_0$, is quite easy.
Let $U$ be a strongly uniform module of dimension
$\aleph_0$. It is easy to see that either

\smallskip

(i) there is a c.d.c\. $\Cal U$ of length
$\omega$ in $U$ such that $U_\alpha$ is complemented over
$U_\beta$ for all $\alpha < \beta < \omega$, or

(ii) there is a c.d.c\. $\Cal U$ of length $\omega$ in $U$ such that
$U_\alpha$ is not complemented over $U_\beta$ for all
$0 \neq \alpha < \beta < \omega$.

\smallskip

In the former case, $U$ is called {\it complementing}; in the
latter, $U$ is {\it narrow}. We refer to \cite{T1, \S2} and
\cite{ET, \S2} for properties and constructions of complementing
and narrow modules of dimension $\aleph_0$.

\medskip

For the more complex case of dimension $d \geq \aleph_1$,
we employ the method of $\Gamma$-invariants as in \cite{T1, \S2}:

\definition{Definition 1.3} Let $\kappa$ be a regular uncountable
cardinal. For any $E \subseteq \kappa$, define
$$\bar E = \{ D \subseteq \kappa \mid \exists C \subseteq \kappa :
\text{ $C$ closed and unbounded in $\kappa$ } \,\,\, \& \,\,\,
D \cap C = E \cap C \}.$$
So $\bar E \in B(\kappa)$.

Let $U$ be a strongly uniform module of
dimension $\kappa$. Let $\Cal U = ( U_\alpha \mid \alpha <
\kappa )$ be a c.d.c\. in $U$. Let
$$E_{\Cal U} = \{ \alpha < \kappa \mid
\exists \beta : \alpha < \beta < \kappa \,\,\, \& \,\,\,
\text{ $U_\alpha$ is not complemented over $U_\beta$}\}.$$
Define $\Gamma (U) = \overline{E_{\Cal U}}$. By \cite{ET, Lemma 1.8},
$\Gamma (U)$ does not depend on the particular choice of the c.d.c\.
$\Cal U$.

$\Gamma (U)$ is called the {\it $\Gamma$-invariant value} of $U$.
We denote by
$\Spec _\Gamma(\kappa,R)$ the {\it $\kappa$-spectrum} of
$\Gamma$, i.e., the set of all $i \in B(\kappa)$ such that there
is a strongly uniform module $U \in \rmod R$ with $\Gamma (U) =
i$. If $\Cal R$ is a class of rings we define
$\Spec _\Gamma(\kappa,\Cal R) = \bigcup_{R \in \Cal R} \Spec
_\Gamma(\kappa,R)$,
the {\it $\kappa$-spectrum} of $\Gamma$ for $\Cal R$.
A $\kappa$-spectrum is said to be {\it full} provided that it is equal
to the whole of $B(\kappa)$.

\enddefinition

The size of $\Spec _\Gamma(\kappa,\Cal R)$ depends substantially
on the properties of $\Cal R$:

\proclaim{Theorem 1.4}
(i) $\Spec _\Gamma(\kappa,\Cal R) = \{ \bar \kappa \}$
for all $\kappa > \aleph_0$ provided that $\Cal R$ is the class of
all commutative rings or $\Cal R$ is the class of all rings with
right Krull dimension.

(ii) For any field $F$, $\Spec _\Gamma(\aleph_1,\Cal R)$ is full
provided that $\Cal R$ is the class of all locally matricial
$F$-algebras.
\endproclaim

\demo{Proof} (i) is by \cite{T1, Theorems 2.10 and 2.12},
and (ii) by \cite{T2, Theorem 2.4}.
\qed\enddemo

The proof of (ii) makes use of a much stronger result, namely
of a realization theorem for ideal lattices of bounded distributive
lattices of size $\leq \aleph_1$ by ideal lattices of von Neumann regular
rings (cf\. \cite{GW, Theorem 4.4} and \cite{T2, Theorem 2.4}).
In particular, the strongly uniform modules are
constructed with the additional property that they are
{\it distributive}, that is, their submodule lattices are distributive.

Nevertheless, by a result of Wehrung \cite{W, Corollary 2.5},
the proof of (ii) does not extend to any $\kappa > \aleph_1$
(see also \cite{PTW, Corollary 4.4}).
It remains open whether $\Spec _\Gamma(\kappa,\Cal R)$ is full for
some $\kappa > \aleph_1$ where $\Cal R$ is the class of all von Neumann
regular rings. \footnote"*"{\it Added in proof: \rm By a different approach,
Pavel R\accent23 u\v zi\v cka recently proved that the spectrum is full
for any regular uncountable cardinal $\kappa$ when $\Cal R$ is the class of
all locally matricial algebras (see also footnote **).}

\bigskip

\heading 2. \, Fullness of the $\kappa$-spectra  \endheading

\medskip

In this section, we will prove that
the $\Spec _\Gamma(\kappa,\Cal R)$ is full for each
regular $\kappa \geq \aleph_1$ where $\Cal R$ is the class
of all rings:

\smallskip

Let $F$ be a field and $\kappa$ be a regular
uncountable cardinal. Fix $S \subseteq \kappa$ with $0 \in S$.

For each $\alpha < \kappa$, put

$$Y_\alpha = \{\langle (\alpha_i,\beta_i); i \leq n \rangle \mid n < \omega;
\alpha_n = \alpha; \alpha_i < \beta_i < \kappa \text{ for all } i \leq
n;$$
$$\alpha_i \in S \text{ for all } 0 < i \leq n; \, \alpha_i < \alpha_{i + 1}
\text{ for all } i < n \}.$$

Observe that $Y_\alpha = \{\langle (\alpha, \beta)\rangle \mid \alpha
< \beta  < \kappa \}$ if $\alpha \notin S$.

For each sequence $y \in Y_\alpha$, $y = \langle (\alpha_i,\beta_i);
i \leq n \rangle$, put $\amax (y) = \alpha_n$, and $\bmax
(y) = \max_{i \leq n} \beta_i \, ( > \amax (y) )$.

Let $Y_{< \alpha} = \cup_{\beta < \alpha} Y_\beta$ and
$Y_{\geq \alpha} = \cup_{\alpha \leq \beta < \kappa} Y_\beta$.
Put $Y = \cup_{\alpha < \kappa} Y_\alpha$. Note that $\card (Y) =
\kappa$.

Denote by $L$ the $F$-linear space with the $F$-basis
$\{ x_\eta \mid \eta \in Y \}$, so

$$L = \bigoplus_{\eta \in Y} F x_\eta$$

has dimension $\kappa$. For each $\alpha < \kappa$,
denote by $L_\alpha$ the $F$-subspace of $L$ generated by
$\{ x_\eta \mid \eta \in Y_{\geq \alpha} \}$. For $\alpha < \beta
< \kappa$ and $\alpha \in S$, we define a subspace

$$L_{\alpha\beta} = \bigoplus_{\eta \in Y_{<\alpha}} F (x_\eta -
x_{\eta\smallfrown(\alpha,\beta)}) \, \oplus \, L_\beta .$$

\definition{Definition 2.1} Let $\nu, \rho \in Y$ be such
that

$$\amax (\rho) \geq \bmax (\nu). \tag{$\ast$}$$

We will define $T = T_{\nu\rho} \in \End_F(L)$.
For $\eta \in Y$, $x_\eta T$ will always be zero or
$x_\theta$, where $\rho$ is an initial segment of $\theta$
which is defined by induction as follows

\smallskip

- if $\nu$ is not an initial segment of $\eta$ then $x_\eta
T = 0$;

- if $\eta = \nu$ then $x_\eta T = x_\rho$;

- if $\nu$ is a proper initial segment of $\eta$, so $\eta =
\eta '\smallfrown (\alpha,\beta)$ and $\nu$ is an initial segment of
$\eta '$, we have $x_{\eta '} T = x_{\rho '}$ for some
$\rho ' \in Y$. If $\rho ' \in Y_{\geq \alpha}$,
we define $x_\eta T = x_{\rho '}$.
If $\rho ' \in Y_{< \alpha}$, we define
$x_\eta T = x_{\rho '\smallfrown (\alpha,\beta)}$.

\smallskip

Denote by $R$ the unital $F$-subalgebra of $\End_F(L)$ generated by
the set $\{ T_{\nu\rho} \mid \nu, \rho \in Y, \amax (\rho)
\geq \bmax (\nu) \}$. Then $L = L_0$ is canonically a (right $R$-) module.
\enddefinition

\proclaim{Lemma 2.2} (i) $L_\alpha$ is a submodule of $L$
for each $\alpha < \kappa$. Moreover, we have $L_\alpha =
x_{\langle (0,1),(\alpha,\alpha + 1)\rangle}R$ \,
for each $0 \neq \alpha \in S$.

(ii) $L_{\alpha\beta}$ is a submodule of $L$
for all $\alpha < \beta < \kappa$ such that $\alpha \in S$.
\endproclaim

\demo{Proof} Let $T = T_{\nu\rho}$, where $\nu, \rho \in
Y$ satisfy ($\ast$).

\smallskip

(i) Let $\eta \in Y_{\geq \alpha}$.
If $\nu \in Y$ is not an initial segment of $\eta$ then $x_\eta
T = 0$. If $\eta = \nu$ then $x_\eta T = x_\rho
\in L_\alpha$ by the assumption ($\ast$).

Let $\nu$ be a proper initial segment of $\eta$, so $\eta =
\eta '\smallfrown (\alpha ',\beta ')$ for some $\alpha \leq
\alpha ' < \beta '$, $\nu$ is an initial segment of $\eta '$, and
$x_{\eta '}T = x_{\rho '}$ for some $\rho ' \in Y$.

If $\rho ' \in Y_{< \alpha'}$ then
$x_\eta T = x_{\rho '\smallfrown (\alpha ',\beta ')} \in L_\alpha$.
If $\rho ' \in Y_{\geq \alpha '}$
then $x_\eta T = x_{\rho '} \in L_\alpha$.

For $0 \neq \alpha \in S$, let $\mu = \langle (0,1)\rangle$ and
$\mu ' = \langle (0,1),(\alpha,\alpha + 1)\rangle$.
Then for each $\eta \in Y_{\geq \alpha + 1}$, we have
$x_\eta = x_{\mu '} T_{\mu ' \eta}$. Similarly, for each
$\eta \in Y_\alpha$ we have $x_\eta = x_{\mu '} T_{\mu \eta}$.

\smallskip

(ii) In view of (i), it suffices to prove that $(x_\eta -
x_{\eta \smallfrown (\alpha,\beta)})T \in L_{\alpha\beta}$ for all
$\eta \in Y_{< \alpha}$.

If $\nu$ is not an initial segment of $\eta \smallfrown (\alpha,\beta)$
then $(x_\eta - x_{\eta \smallfrown (\alpha,\beta)})T = 0$.

If $\nu = \eta \smallfrown (\alpha,\beta)$ then $(x_\eta - x_{\eta \smallfrown
(\alpha,\beta)})T
= (- x_{\eta \smallfrown (\alpha,\beta)})T = - x_\rho \in L_\beta$
by ($\ast$).

If $\nu = \eta$ then $(x_\eta - x_{\eta \smallfrown (\alpha,\beta)})T
= x_\rho - (x_{\eta \smallfrown (\alpha,\beta)})T$. If $\rho \in Y_{< \alpha}$,
then $(x_{\eta \smallfrown (\alpha,\beta)})T = x_{\rho \smallfrown
(\alpha,\beta)}$,
so $(x_\eta - x_{\eta \smallfrown (\alpha,\beta)})T \in L_{\alpha\beta}$.
If $\rho \in Y_{\geq \alpha}$, then
$(x_{\eta \smallfrown (\alpha,\beta)})T = x_\rho$, so
$(x_\eta - x_{\eta \smallfrown (\alpha,\beta)})T = 0$.

Assume that $\nu$ is a proper initial segment of $\eta$,
so $\eta = \eta '\smallfrown (\alpha ',\beta ')$ for some $\alpha ' <
\alpha$, and $\nu$ is an initial segment of $\eta '$.
We have $x_{\eta '}T = x_{\rho '}$ where
$\rho ' \in Y$.

If $\rho ' \in Y_{< \alpha '}$ then
$x_\eta T = x_{\rho '\smallfrown (\alpha ',\beta')}$ while
$(x_{\eta \smallfrown (\alpha,\beta)})T =
x_{\rho '\smallfrown (\alpha ',\beta') \smallfrown (\alpha,\beta)}$,
because $\alpha' < \alpha$. So
$(x_\eta - x_{\eta \smallfrown (\alpha,\beta)})T \in L_{\alpha\beta}$.

Assume $\rho ' \in Y_{\geq \alpha '}$,
so $x_\eta T = x_{\rho '}$.
If $\rho ' \in Y_{< \alpha}$, then
$(x_{\eta \smallfrown (\alpha,\beta)})T = x_{\rho '\smallfrown
(\alpha,\beta)}$,
so $(x_\eta - x_{\eta \smallfrown (\alpha,\beta)})T \in L_{\alpha\beta}$.
If $\rho ' \in Y_{\geq \alpha}$, then
$(x_{\eta \smallfrown (\alpha,\beta)})T = x_{\rho '}$,
so $(x_\eta - x_{\eta \smallfrown (\alpha,\beta)})T = 0$.
\qed\enddemo

\proclaim{Lemma 2.3} $\Cal L = ( L_\alpha \mid \alpha < \kappa )$
is a c.d.c\. in $L$. \endproclaim

\demo{Proof} Clearly, $\Cal L$ is strictly decreasing and
continuous. Let $X$ be a non-zero submodule of $L$ and take $0
\neq x \in X$. So $x = \sum_{\eta \in Y} f_\eta x_\eta$ and
$f_\eta = 0$ for almost all, but not all, $\eta \in Y$.
Take $\nu \in Y$ such that $f_\nu \neq 0$ and $\nu$ is not a
proper initial segment of any $\eta \in Y$ with $f_\eta \neq 0$.
Let $\alpha = \bmax (\nu)$.
Take any $\rho \in Y_{\geq \alpha}$ and let $T = T_{\nu\rho}$.
Then $x T = (f_\nu x_\nu) T = f_\nu x_\rho$, so $x_\rho \in X$.
This proves that $L_\alpha \subseteq X$, and $\Cal L$ is
cofinal.
\qed\enddemo

\proclaim{Proposition 2.4} Let $\gamma < \kappa$. Then
$\Cal L _\gamma = ( L_\alpha \mid \gamma \leq \alpha < \kappa )$
is a c.d.c\. in $L_\gamma$ such that $E_{\Cal L _\gamma} =
[\gamma,\kappa) \setminus S$. In particular, $\Gamma (L_\gamma) =
\overline{\kappa \setminus S}$.
\endproclaim

\demo{Proof} By Lemma 2.3, $\Cal L _\gamma$ is a c.d.c\. in
$L_\gamma$.

We prove that $L_\alpha$ is complemented
over $L_\beta$ in $L_\gamma$ provided that
$\gamma < \alpha < \beta < \kappa$ and $\alpha
\in S$. By modularity, it is enough to prove this for $\gamma =
0$:

Clearly, $L = L_\alpha + L_{\alpha\beta}$.
Take $x \in L_\alpha \cap L_{\alpha\beta}$. Then $x = y +
z$, where $y \in \bigoplus_{\eta \in Y_{< \alpha}} F (x_\eta -
x_{\eta \smallfrown (\alpha,\beta)})$ and $z \in L_\beta$. Since
$x \in L_\alpha$, we have $y = 0$, so
$L_\beta = L_\alpha \cap L_{\alpha\beta}$.

It remains to prove that $L_\alpha$ is not complemented over
$L_\beta$ in $L_\gamma$ provided that $\gamma < \alpha < \beta < \kappa$
and $\alpha \notin S$:

Assume there is a submodule $X$ in $L$ such that
$L_\gamma = L_\alpha + X$ and $L_\beta = L_\alpha \cap X$.
Let $\nu = \langle (\gamma,\gamma + 1)\rangle \in Y_{\geq \gamma}$,
$\rho = \langle (\alpha,\alpha + 1)\rangle \in Y_{\geq \alpha}$ and
take $T = T_{\nu\rho}$.
By assumption,
there are $x \in X$ and $y \in L_\alpha$ such that $x_\nu = x + y$.
Since $\alpha \notin S$, we have
$(L_\alpha) T \subseteq L_{\alpha + 1}$. So
$x T = x_\rho - y T \in L_\alpha \setminus L_{\alpha + 1}$.
On the other hand, $x T \in X$, so $x T \in L_\beta$, a contradiction.
\qed\enddemo

The following lemma says that each $L_\alpha$, $\alpha < \kappa$,
is a rigid module in the sense that $\End _R(L_\alpha)$ is minimal possible.

\proclaim{Lemma 2.5} $\End _R(L_\alpha) = F$ for all $\alpha <
\kappa$. \endproclaim

\demo{Proof} Let $0 \neq e \in \End_R(L_\alpha)$.

First, we prove that $\ker e = 0$. If not, by Lemma 2.3,
there is $\beta < \kappa$ such that $L_\beta
\subseteq \ker e \cap \Image e$. Take $\nu \in Y_{\geq \beta}$.
Let $x \in L_\alpha$ be such that $e x = x_\nu$. Then
$x = \sum_{\eta \in Y_{\geq \alpha}} f_\eta x_\eta$, and the set
$A = \{ \eta \in Y_{\geq \alpha} \mid f_\eta \neq 0 \}$ is finite.
W.l.o.g., we may assume that $e x_\eta \neq 0$ for all $\eta \in
A$. Then, for each $\eta \in A$, $\nu$ is not an initial segment
of $\eta$. Take $\rho \in Y_{\geq \beta}$ such that ($\ast$) holds.
Put $T = T_{\nu\rho}$. Then $0 = e(xT) = (ex)T = x_\rho$, a
contradiction.

Next, we prove that for each $\eta \in Y_{\geq \alpha}$, there is $f_\eta \in
F$ such that $e x_\eta = f_\eta x_\eta$. Clearly, $e x_\eta =
\sum_{\tau \in Y_{\geq \alpha}} f_\tau x_\tau$, and the set
$A = \{ \tau \in Y_{\geq \alpha} \mid f_\tau \neq 0 \}$ is finite.
Since $\ker e = 0$, at least one $\tau \in A$ must contain
$\eta$ as an initial segment. Let $\tau _0 \in A$ be maximal such.
If $\tau _0 \neq \eta$, then taking $\rho \in Y_{\geq \alpha}$ such that
$\amax(\rho) \geq \bmax(\tau _0)$, we see that $T_{\tau_0 \rho}$ maps
$e x_\eta$ to $f_{\tau_0} x_\rho$, while $x_\eta T_{\tau_0 \rho} = 0$,
a contradiction. This shows that $\tau _0 = \eta$.

Let $\tau \in A \setminus
\{ \eta \}$ be maximal. If $\tau$ is not an initial segment of
$\eta$, then taking $\rho \in Y_{\geq \alpha}$ such that
$\amax(\rho) \geq \bmax(\tau)$, we see that
$T_{\tau\rho}$ maps $x_\eta$ to $0$, but
$(e x_\eta) T_{\tau\rho} = f_\tau x_\rho$, a contradiction.
So $\tau$ is a proper initial segment of $\eta$, $\eta = \eta '
\smallfrown (\beta,\gamma)$, and $\tau$ is an initial segment of $\eta '$.
Take $\rho \in Y$ such that $\amax (\rho) \geq \bmax
(\eta)$ and let $T = T_{\tau\rho}$.
Then $x_{\eta '}T = x_{\rho '}$ for some
$\rho '$ containing $\rho$ as an initial segment. Then $x_\eta T
= x_{\rho '}$, so $e x_{\rho '} = (e x_\eta) T = f_\tau x_\rho +
f_\eta x_{\rho '}$. On the other hand, $e x_{\rho '} = (e x_\eta)
T_{\eta\rho '} = f_\eta x_{\rho '}$. So $f_\tau = 0$, a
contradiction.

Finally, we prove that $f_\nu = f_\rho$ for all $\nu, \rho \in
Y_{\geq \alpha}$.
This is clear when ($\ast$) holds. But then $f_\nu = f_\rho = f_\nu '$,
where $\nu, \nu ' \in Y_{\geq \alpha}$ are arbitrary, and $\rho =
\langle (\beta,\gamma)\rangle $ is such that ($\ast$) holds and
$\beta = \amax (\rho) \geq \bmax (\nu ')$.
\qed\enddemo

\proclaim{Theorem 2.6} Let $\kappa$ be a regular uncountable
cardinal and $i \in B(\kappa)$. Let $F$ be a field and $L$ be an
$F$-linear space of dimension $\kappa$.

Then there exists an $F$-subalgebra, $R$, of $End_F(L)$ such that
$L$, viewed as a right $R$-module, is strongly uniform with
$\Gamma (L) = i$ and $\End_R(L) = F$. \endproclaim

\demo{Proof} By Proposition 2.4 and Lemma
2.5. \qed\enddemo

In the construction of Theorem 2.6, different elements of
$B(\kappa)$ occur as values of the $\Gamma$-invariant of modules over
different algebras. This is easily improved in our main result

\proclaim{Theorem 2.7} Let $\lambda$ be an uncountable
cardinal and $F$ be a field. Then there exists an $F$-algebra $R$
such that for any regular uncountable cardinal $\kappa \leq
\lambda$ and any $i \in B(\kappa)$ there is a strongly uniform
module $L \in \rmod R$ such that $\End _R (L) = F$ and
$\Gamma (L) = i$. In particular, $\Spec _\Gamma(\kappa,\Cal R)$ is
full. \endproclaim

\demo{Proof} For each regular uncountable cardinal $\kappa
\leq \lambda$ and each $i \in B(\kappa)$, denote by $R_{\kappa i}$ the
$F$-algebra, and by $L_{\kappa i}$ the right $R_{\kappa
i}$-module, constructed in Theorem 2.6. Let $R = \prod_{\kappa,
i} R_{\kappa i}$ (the ring direct product).
Then each $L = L_{\kappa i}$ is canonically a right
$R$-module, and the $R$- and $R_{\kappa i}-$ submodule lattices
of $L$ coincide. It follows that $\Gamma (L) = i$. Moreover,
$\End _R(L) = \End _{R_{\kappa i}} (L) = F$.
\qed\enddemo

\bigskip

\heading 3. \, The $\Gamma$-invariant of two-sided ideal lattices  \endheading

\medskip

The $\Gamma$-invariant of strongly uniform modules as defined in Section 1
is completely determined by properties of submodule
lattices of the respective modules. In fact, this is a particular instance
of a more general $\Gamma$-invariant, the $\Gamma$-invariant of
strongly dense lattices \cite{ET, \S1}.

Recall that a bounded modular lattice
$( A, \wedge, \vee, 0,  1)$ is {\it strongly dense}
provided that it contains a continuous strictly decreasing cofinal
chain ({\it c.d.c.}) consisting of non-zero elements of $A$.
If $0 \neq b < a \leq 1 \in A$, then
$a$ is {\it complemented} over $b$ provided that there exists $c \in
A$ with $a \wedge c = b$ and $a \vee c = 1$. As in Definition 1.3,
we can define for each c.d.c. $\Cal U$ of length $\kappa$ in $A$ the set
$\overline{E_{\Cal U}} \in B(\kappa)$.
Then $\Gamma (A) = \overline{E_{\Cal U}}$ does not depend on the choice of
the c.d.c. $\Cal U$, and it is called the {\it $\Gamma$-invariant value} of
the lattice $A$, \cite{ET, \S1}.

This $\Gamma$-invariant is of particular interest in the case
when $A = L_2(S)$, the two-sided ideal lattice of an algebra $S$.
Indeed, the proof of Theorem 1.4(ii) makes essential use of this
case: for $\kappa = \aleph_1$, applying a construction due to
Goodearl and Wehrung \cite{GW, Theorem 4.4} together with
\cite{ET, Theorem 1.15}, one can realize each $i \in
B(\aleph_1)$ as $\Gamma (L_2(S))$ for a locally matricial
$F$-algebra $S$. In particular, $L_2(S)$ is a distributive
lattice. Let $R = S \otimes _F S^{\text{op}}$, where
$S^{\text{op}}$ is the opposite $F$-algebra of $S$. Then $S$ is
a (right $R$-) module whose submodule lattice is canonically isomorphic
to $L_2(S)$. So $S$ is a strongly uniform module
of dimension $\kappa$. Moreover, $i = \Gamma(L_2(S)) = \Gamma
(S)$, so $i$ is realized as the $\Gamma$-invariant value of a
distributive strongly uniform module.

For $\kappa > \aleph_1$, the question of the possible
values of the $\Gamma$-invariant of strongly dense two-sided
ideal lattices remains open. \footnote"**"{\it Added in proof: \rm
Recently, Pavel R\accent23 u\v zi\v cka proved that the ideal
lattice of any bounded distributive lattice is isomorphic to the
lattice of two-sided ideals of a locally matricial algebra.
>From \cite{ET, Theorem 1.15}, it easily follows that
the spectrum of the $\Gamma$-invariant of strongly dense two-sided
ideal lattices is full for any $\kappa > \aleph_1$.
More details appear in R\accent23 u\v zi\v cka's manuscript
"Lattices of two-sided ideals of locally matricial algebras
and the $\Gamma$-invariant problem".} Nevertheless, Theorem 2.6
provides a realization of any $i \in B(\kappa)$ as $\Gamma(A)$ where
$A$ is a lower interval in $L_2(S)$ for an $F$-algebra $S$:

\proclaim{Corollary 3.1} Let $F$ be a field, $\kappa$ be a
regular uncountable cardinal, $i \in B(\kappa)$,
$R$ be the $F$-algebra and $L$ be the module constructed
in Theorem 2.6. Let

$$S = \{ \left ( \smallmatrix
f \, l \\ 0 \, r \endsmallmatrix \right ) \mid f \in F, l \in L,
r \in R \}.$$

Let $I = \{ \left ( \smallmatrix
0 \, l \\ 0 \, 0 \endsmallmatrix \right ) \mid l \in L \}$. Then
$S$ is an $F$-algebra and $I \in L_2(S)$.
Denote by $A$ the interval in $L_2(S)$ consisting
of all two-sided ideals contained in $I$.
Then $A$ is a strongly dense lattice of dimension $\kappa$
and $\Gamma (A) = i$.
\endproclaim

\demo{Proof} Clearly, $A$ is isomorphic to the (right $R$-) submodule
lattice of $L$, so the assertion follows by Theorem 2.6. \qed\enddemo

Though our construction in Section 2 applies to an arbitrary
regular uncountable cardinal $\kappa$, it neither produces $R$ which is
von Neumann regular nor $L$ which has a distributive lattice of
submodules. So Theorem 1.4(ii) provides a stronger result
in the particular case of $\kappa = \aleph_1$:

\proclaim{Lemma 3.2} Neither of the algebras $R$ appearing
in Theorems 2.6 and 2.7 is von Neumann regular. Neither of the
strongly uniform modules $L$ from Theorems 2.6 and 2.7
is distributive. \endproclaim

\demo{Proof} To see that $R$ in Theorem 2.6 (and hence in 2.7)
is not von Neumann regular take $\alpha + 1 <
\beta < \kappa$, $\mu = \langle (\alpha,\alpha + 1)\rangle$ and
$\phi = \langle (\alpha + 1, \beta)\rangle$.
Then $T_{\mu\phi}$ has no pseudo-inverse in $R$.

Indeed, if $T \in R$ is such that $T_{\mu\phi} T T_{\mu\phi} =
T_{\mu\phi}$, then $x_\phi T \in T^{-1}_{\mu\phi}(x_\phi) \cap
L_{\alpha + 1}$ by Lemma 2.2(i). It follows that $x_\phi T = x_\tau$,
where $\tau = \mu \smallfrown (\alpha + 1,\gamma)$ for some
$\alpha + 1 < \gamma < \kappa$.
Now, any $T_{\nu\rho}$, with $\nu, \rho \in Y$ satisfying
($\ast$), maps $x_\phi$ to zero or to $x_{\phi '} \in L_{\alpha +
2}$. On the other hand, we have $T = f.1 + t \in R$, where $f \in F$ and
$t$ is an $F$-linear combination of finite products of elements of the form
$T_{\nu\rho}$, with $\nu, \rho \in Y$ satisfying ($\ast$).
Then $x_\tau = x_\phi T = fx_\phi + x_\phi t$, where $x_\phi t \in
L_{\alpha + 2}$, a contradiction.

To see that the module $L$ in Theorem 2.6 (and hence in 2.7)
is not distributive, fix
$\alpha < \kappa$, and for each $\alpha + 1 < \beta < \kappa$
let $\phi _\beta = \langle (\alpha + 1, \beta)\rangle$.
Then $(x_{\phi _\beta} + L_{\alpha + 2}) r = f x_{\phi _\beta} +
L_{\alpha + 2}$ for any $r = f.1 + t \in R$, where $f \in F$ and $t$ is
an $F$-linear combination of finite products of elements of the form
$T_{\nu\rho}$, with $\nu, \rho \in Y$ satisfying ($\ast$).
So the $R$-submodules, and the $F$-subspaces, of
$N_\alpha = \bigoplus_{\alpha + 1 < \beta < \kappa} (x_{\phi
_\beta} + L_{\alpha + 2})R \subseteq L/L_{\alpha + 2}$ coincide.
Since $\dim _F(N_\alpha) = \kappa > 1$, the module $N _\alpha$,
and hence $L$, is not distributive.
\qed\enddemo

The results above suggest the question of the structure of $L_2(R)$
for the $F$-algebra $R$ constructed in Theorem 2.6. We will prove
that $L_2(R)$ is strongly dense, but in contrast with
the Goodearl-Wehrung construction, $L_2(R)$ is always narrow.
First, we need more information about the arithmetic of the
algebra $R$:

\smallskip

Let $\nu, \nu', \rho, \rho' \in Y$ be such
that ($\ast$) holds and $\amax (\rho ') \geq \bmax (\nu')$.
We will compute $T_{\nu \rho} T_{\nu ' \rho '}$:

(I) If $\nu'$ is not an initial segment of $\rho$ and $\rho$ is not
an initial segment of $\nu '$, then $T_{\nu \rho} T_{\nu ' \rho
'} = 0$.

(II) If $\rho = \nu ' \smallfrown \tau$, then $T_{\nu \rho}
T_{\nu ' \rho '} = T_{\nu , \rho ' \smallfrown \tau '}$ where
$\tau ' = \emptyset$ provided that $\amax (\rho) \leq \amax (\rho ')$,
and $\tau '$ is the final segment of $\tau$ consisting of all pairs
whose first component is $> \amax(\rho ')$ provided that
$\amax (\rho) > \amax (\rho ')$.

(III) If $\nu ' = \rho \smallfrown \tau$ and $\tau \neq \emptyset$,
then

$$T_{\nu \rho} T_{\nu ' \rho '} = \bigoplus _{\sigma \in X}
T_{\nu \smallfrown \sigma \smallfrown \tau , \rho '},$$

where $X$ consists of the empty set and of all
elements of $Y_{\leq \amax(\rho)}$ whose initial pair has first
component $> \amax (\nu)$.

Further, let $t = \prod _{i \leq n} T_{\nu_i \rho _i}$ where
$n < \omega$ and $\amax (\rho _i) \geq \bmax (\nu _i)$ for all $i
\leq n$. If $n > 0$ and $t$ is irredundant (in the sense that the product
cannot be simplified using (II) for successive factors),
then (III) shows that

$$t = \bigoplus _{\sigma_0 \in X_0, \dots, \sigma _n \in X_n}
T_{\nu_0 \smallfrown \sigma_0 \smallfrown \tau_0 \smallfrown \dots
\smallfrown \sigma_n \smallfrown \tau_n  , \rho _n} $$

where $\nu_{i + 1} = \rho_i \smallfrown \tau_i$ for all $i < n$,
$\tau_i \neq \emptyset$ for all $i \leq n$, and for each $i
\leq n$, $X_i$ consists of the empty set and of all elements of
$Y_{\leq \amax(\rho_i)}$ whose initial pair has first component
$> \amax (\nu_i)$.
Note that $\amax(\nu _0) < \amax(\rho _0) < \amax (\nu_1)
\dots < \amax (\rho_{n -1}) < \amax (\nu_n) < \amax (\rho_n)$.

Let $r \in R$. Then $r$ can be expressed as an $F$-linear
combination

$$r = f.1 + \sum_{j < m} f_j t_j,  \tag{$\ast \ast$}$$

where $m < \omega$, $f \in F$, $0 \neq f_j \in F$ and $t_j$
is a finite irredundant product of elements of the form
$T_{\nu\rho}$ with $\nu, \rho \in Y$ satisfying ($\ast$)
for each $j < m$.

So each $t_j$ is of the form

$$t_j = \bigoplus _{\sigma_{j0}\in X_{j0}, \dots, \sigma
_{jn_j}\in X_{n_j}} T_{\nu_{j0} \smallfrown \sigma_{j0}
\smallfrown \tau_{j0} \smallfrown
\dots \smallfrown \sigma_{jn_j} \smallfrown \tau_{jn_j} ,
\rho_{jn_j}}$$

(in order to unify our notation, we set $n_j = 0$, $X_j = \{
\emptyset \}$ and $\tau_j = \emptyset$ in the case when
$t_j = T_{\nu_{j0} \rho_{j0}}$ has exactly one factor).

We will say that ($\ast \ast$) is a {\it canonical form}
of $r$ provided that each $t_j$ is irredundant and $t_j \neq
t_{j'}$ for all $j \neq j' < m$.

\smallskip

\proclaim{Theorem 3.3} $L_2(R)$ is a strongly dense lattice of
dimension $\kappa$ and $\Gamma(L_2(R)) = \bar \kappa$.
\endproclaim

\demo{Proof} For each $\alpha < \kappa$, define

$$I_\alpha = \{ r \in R \mid \Image r \subseteq L_\alpha \}.$$

The proof is divided into three lemmas:
\enddemo

\proclaim{Lemma 3.4} Let $r \in R$ be in the canonical
form ($\ast \ast$). Let $\alpha > 0$. Then $r \in I_\alpha$ iff $f =
0$ and $\amax (\rho_{jn_j}) \geq \alpha$ for all $j < m$.
In particular, $I_\alpha$ coincides with the ideal of $R$ generated
by all $T_{\nu\rho}$ such that $\nu, \rho \in Y$ satisfy
($\ast$) and $\amax(\rho) \geq \alpha$.
\endproclaim

\demo{Proof} The `if' part is clear, since $r$ then
maps into $L_\alpha$.

For the `only if' part, assume that $r \in L_\alpha$.
If $f \neq 0$ then we take $\eta \in Y_0$ such that $\nu_{j0}$
is not an initial segment of $\eta$ for all $j < m$.
Then $x_\eta r = f x_\eta \notin L_\alpha$, a contradiction.

Proving indirectly, we can w.l.o.g\. assume that $f = 0$ and
$\amax (\rho _{jn_j}) < \alpha$ for all $j < m$.
Let $i < m$ be such that $\rho _{in_i}$ is minimal.
Since $r \in I_\alpha$, we have $\card(J) \geq 2$ where
$J = \{ j < m \mid \rho _{jn_j} = \rho _{in_i} \}$.
Let $j \in J$ be such that
$\nu _{j0}$ is minimal. Since $r \in I_\alpha$, we have
$\card(J_0) \geq 2$ where
$J_0 = \{ j' \in J \mid \nu _{j'0} = \nu _{j0} \}$.

If there is $k \in J_0$ such that $n_{k0} = 0$, then
there exists $k' \in J_0$ such that $k' \neq k$ and $t_{k'} = t_k
= T_{\nu_{k0} \rho_{k0}}$
which contradicts the assumption that $(\ast \ast)$ is canonical.

Otherwise, let $k \in J_0$ be such that
$\amax(\sigma_{k0})$ is maximal. Then $\card(J_1) \geq 2$ where
$J_1 = \{ k' \in J_0 \mid \amax(\rho_{k'0}) = \amax(\rho_{k0})
\}$. Let $l \in J_1$ be such that
$\amax(\tau_{l0})$ is minimal. Then $\card (J_2) \geq 2$ where
$J_2 = \{ l' \in J_1 \mid \tau _{l'0} = \tau _{l0} \}$. Proceeding
similarly, after finitely many steps we obtain a pair $j \neq j'
< m$ such that $t_j = t_{j'}$ which contradicts the assumption
that $(\ast \ast)$ is canonical.
\qed\enddemo

Note that Lemma 3.4 implies that the canonical form
($\ast \ast$) is unique for each $r \in R$. That is, all the
irredundant products together with $1 \in R$ form an $F$-basis of $R$.

\proclaim{Lemma 3.5} $\Cal I = (I_\alpha \mid \alpha < \kappa)$ is a
c.d.c\. in $L_2(R)$.
\endproclaim

\demo{Proof} Clearly, $I_\alpha \in L_2(R)$. Since
$T_{\langle (0,1)\rangle , \langle (\alpha,\alpha +1)\rangle}
\in I_\alpha \setminus I_{\alpha
+ 1}$, $\Cal I$ is strictly decreasing. By definition, $I_\alpha
= \cap_{\beta < \alpha} I_\beta$ for all limit ordinals $\alpha <
\kappa$, so $\Cal I$ is continuous.

Let $0 \neq r \in R$. We will prove that $T_{\nu\rho} \in RrR$
for some $\nu, \rho \in Y$ satisfying ($\ast$).
Consider the canonical form of $r$, ($\ast \ast$).

If $f \neq 0$ then there is $T_{\nu\rho}$ satisfying ($\ast$)
such that $\nu \in Y_0$ is not an initial segment of
$\rho _{jn_j}$ for any $j < m$.
Then $rT_{\nu\rho} = fT_{\nu\rho}$, so $T_{\nu\rho} \in RrR$.

Assume $f = 0$. Multiplying $r$ by an appropriate $T_{\nu' \rho'}$
on the right and using (III), we can w.l.o.g. assume that
$\rho ' = \rho_{jn_j}$ for all $j < m$. Since ($\ast \ast$) is
canonical, an argument similar to the one in the proof of Lemma
3.4 shows that there exist $\rho '' \in Y$ and $\beta <
\kappa$ such that there is $j < m$ with
$T_{\langle (0,\beta)\rangle \rho ''} r
= T_{\langle (0,\beta)\rangle \rho ''} t_j = T_{\langle (0,\beta)\rangle
\rho}$ where $\rho'$ is an initial segment of $\rho$. Then $T_{\langle
(0,\beta)\rangle \rho} \in RrR$.

Let $s = T_{\nu\rho} \in R$ with $\nu, \rho
\in Y$ satisfying ($\ast$). Put $\alpha = \bmax (\rho)$.
To finish the proof it suffices to show that $I_\alpha \subseteq
RsR$. By Lemma 3.4, it is enough to show that $T_{\nu' \rho '}
\in RsR$ whenever $\nu', \rho' \in Y$ satisfy $\amax
(\rho ') \geq \bmax(\nu ')$ and $\amax(\rho ') \geq \alpha$.

If $\nu \in Y_{\geq 1}$, then $T_{\langle (0,1)\rangle \nu} \, s =
T_{\langle (0,1)\rangle \rho}$ and $T_{\langle (0,1)\rangle \rho}
T_{\rho \rho'} =
T_{\langle (0,1)\rangle \rho'}$, so $T_{\nu' \rho '} = T_{\nu' \rho ''}
T_{\langle (0,1)\rangle ,\rho'} \in RsR$, where $\rho ''$ is obtained from
$\rho '$ by adding (replacing by) the initial pair $(0,1)$.

Let $\nu \in Y_0$ so $\nu = \langle (0, \beta)\rangle$ where
$0 < \beta < \alpha$. As above, we get $T_{\nu\rho'} \in sR$,
and $T_{\nu ' \rho '} \in RsR$.
\qed\enddemo

\proclaim{Lemma 3.6} $\Gamma (L_2(R)) = \bar \kappa$.
\endproclaim

\demo{Proof} Let $0 < \alpha < \beta < \kappa$. Assume there
exists $C \in L_2(R)$ such that $I_\alpha + C = R$ and $I_\alpha
\cap C = I_\beta$. In particular, there exists $r \in I_\alpha$
such that $1 - r \in C$ and $I_\alpha (1 - r) \subseteq I_\alpha
\cap C = I_\beta$. Consider the canonical form of $r$, ($\ast
\ast$). By Lemma 3.4, $f = 0$. Moreover, there exists
$\alpha < \gamma < \kappa$ such that for each $j < m$,
if a pair $(\alpha,\delta_j)$ occurs in $\nu_{j0}$,
then $\gamma \neq \delta_j$. Then $s = T_{\langle (0,1)\rangle , \langle
(\alpha,\gamma)\rangle}
\in I_\alpha \setminus I_\beta$ and $sr = 0$, so $s (1 - r) = s
\in I_\beta$, a contradiction. This proves that $I_\alpha$ is not
complemented over $I_\beta$. By Lemma 3.5, $\Gamma (L_2(R)) =
\bar \kappa$.
\qed\enddemo

\bigskip

\heading 4. \, Complexity of Ziegler spectra of infinite
dimensional algebras \endheading

\medskip
Theorem 2.7 cannot be improved to produce
a proper class of strongly uniform modules with different
values of the $\Gamma$-invariant over a fixed ring $R$:

\proclaim{Lemma 4.1} Let $R$ be a ring. For each right ideal $I$
of $R$ such that $R/I$ is strongly uniform, denote by $d_I$ the
dimension of $R/I$ (for example, $d_I = 1$ for any maximal right
ideal $I$). Let $\kappa _R = \sup_{I} d_I$. Then each strongly uniform
module has dimension $\leq \kappa _R$. \endproclaim

\demo{Proof} Let $U$ be a strongly uniform module of dimension
$\lambda$. Let $E$ be the injective hull of $U$.
Then $E$ is strongly uniform, and has dimension $\lambda$.
On the other hand, $E$ is the injective hull of some cyclic module
$R/I$. Then also $R/I$ is strongly uniform of dimension
$\lambda$, so $\lambda = d_I$.
\qed\enddemo

We do not know whether we can improve Theorem 2.7 to
produce injective uniform (= indecomposable injective) modules with
prescribed values of the $\Gamma$-invariant.
Nevertheless, slightly modifying the invariant, we can produce
the relevant examples:

\definition{Definition 4.2} Let $\kappa$ be a regular uncountable
cardinal. Let $U$ be a strongly uniform module of dimension $\kappa$.
By modularity of submodule lattices, we have $\Gamma (V) \leq \Gamma
(W)$ for any non-zero submodules $V \subseteq W \subseteq U$. So the set
$$\Cal G (U) = \{ \Gamma (V) \mid 0 \neq V \subseteq U \}$$
is a lower directed subset of $B(\kappa)$.

If $U$ is such that $\Cal G (U)$ has a least element, we define

$$\Gamma ^* (U) = \min \Cal G (U);$$

otherwise, $\Gamma ^* (U)$ is not defined.
\enddefinition

Recall that for a ring $R$, the {\it Ziegler spectrum} of $R$,
$\Zg (R)$, is a topological space whose points are
(isomorphism classes of) indecomposable pure-injective modules
and the topology has the property that closed subsets
correspond bijectively to complete theories of modules closed
under products. (A closed subset $C$ corresponds to the complete
theory of the module $M = \oplus_{N \in C} N^{(\omega)}$).
Despite being a set, the Ziegler spectrum captures most model
theoretic properties of the class $\rmod R$, cf\. \cite{P1}, \cite{P2}.

We finish by showing that the point structure of $\Zg (R)$ is
very complex in case $R$ is the infinite dimensional algebra
constructed in Theorem 2.7:

\proclaim{Theorem 4.3} Let $\lambda$ be an uncountable
cardinal and $F$ be a field. Then there exists an $F$-algebra $R$
such that for any regular uncountable cardinal $\kappa \leq
\lambda$ and any $i \in B(\kappa)$ there is a
strongly uniform module $I \in \Zg (R)$ such
that $\Gamma ^* (I) = i$.
\endproclaim

\demo{Proof} Let $R$ be as in Theorem 2.7, and $L = L_{\kappa i} \in \rmod R
\cap \rmod R_{\kappa i}$ be the strongly uniform module constructed in
Theorem 2.7, with $\Gamma (L) = i$. By Proposition 2.4,
the right $R$-submodule $L_\alpha$ has
$\Gamma$-invariant value equal to $i$ for all $\alpha < \kappa$.
Denote by $E$ the injective hull of the right $R$-module $L$.
>From Lemma 2.3 we get that $\{ L_\alpha \mid \alpha
< \kappa \}$ is cofinal in $E$. So $i$ is the least element of
the lower directed set $\Cal G (E)$. This proves that
$\Gamma ^* (E) = \Gamma (L) = i$. Finally, there is a (unique)
element $I \in \Zg (R)$ which is isomorphic to $E$, so
$\Gamma ^* (I) = i$.
\qed\enddemo

\enddocument
\end